\documentclass[a4paper,12pt]{article}
\usepackage{amsmath,amssymb,amsthm,geometry,mathrsfs}
\usepackage{hyperref}
\usepackage[utf8]{inputenc}
\usepackage{listings}
\usepackage{graphicx}
\usepackage{xcolor}
\usepackage{algorithm}
\usepackage{algpseudocode}
\usepackage{enumitem}

\newtheorem{theorem}{Theorem}[section]

\title{Cohomologically calibrated affine connections and the Einstein condition on \( S^2 \times T^2 \)}
\author{Alexander Pigazzini, Magdalena Toda}
\date{August 2025}

\begin{document}

\maketitle

\begin{abstract}

\noindent This paper applies the recently developed framework of cohomologically calibrated affine connections to the fundamental problem of constructing non-Riemannian Einstein manifolds. In this framework, the torsion of a connection is intrinsically related to the global topology of the manifold, represented by the de Rham cohomology class specified by a set of real parameters. We focus on the product manifold \(S^2 \times T^2\), whose third cohomology group is \(H^3(S^2 \times T^2; \mathbb{R}) \cong \mathbb{R}^2\). We analyze how the geometry is modeled by the choice of the torsion tensor \(T\) within the family \(\mathcal{T}_\omega\), defined by the property that each member of this family must have an associated 3-form \(T^\flat\) such that it represents the nontrivial cohomology class via Hodge decomposition. Our analysis reveals a dependence on this choice. First, we show that using a torsion tensor, which produces a strictly positive biorthogonal curvature, leads to a non-diagonal Ricci tensor, creating a structural obstacle to any Einstein solution. Conversely, we then show that using the torsion tensor associated with the purelly harmonic 3-form allow us the construction of an explicit non-Riemannian Einstein solution. Our work thus demonstrates that the cohomologically calibrated affine connections allow a family of feasible connections rich enough to allow for several geometries intrinsically justified by the differential topology of the manifold.
\end{abstract}

\bigskip

\noindent \textbf{Keywords:} Cohomologically calibrated affine connections, Non-Riemannian Einstein manifolds, torsion, cohomology, 4-manifolds, non-simply connected manifold.
\\
\\
\textbf{MSC 2020:} 53C05, 53C15, 53C25, 58A14.

\section{Introduction}

Affine connections with torsion arise naturally in a variety of mathematical and physical contexts, from Cartan geometry \cite{Cartan}, and special holonomy \cite{Kobayashi}, \cite{Chiossi} (the latter directly connects the G-structures to the intrinsic torsion), to string theory \cite{Duff}, \cite{Friedrich}, and generalized geometry \cite{Hitchin}, \cite{Gualtieri}. While the Levi-Civita connection, characterized by being both torsion-free and metric-compatible, remains central in Riemannian geometry, relaxing the torsion-free condition opens the door to a broader class of geometric structures with rich curvature behavior and topological implications \cite{Agricola}. Among these, connections with totally antisymmetric torsion have received particular attention due to their analytical tractability and geometric relevance (see, \cite{Gauduchon}).

In this paper, we investigate whether a specific affine connection with torsion, referred to as the \emph{cohomologically calibrated affine connection} and denoted \( \nabla^{\mathcal{C}} \), can satisfy the Einstein condition on the product manifold \( S^2 \times T^2 \). The connection \( \nabla^{\mathcal{C}} \) was recently introduced in \cite{Pigazzini2025}, where it was constructed on \( S^2 \times T^2 \) by defining a totally antisymmetric torsion tensor \( T \). The associated 3-form \( T^\flat \), defined via the metric \( g \) as \( T^\flat(X,Y,Z) = g(T(X,Y), Z) \), is associed with a non-trivial cohomology class in \( H^3(S^2 \times T^2; \mathbb{R}) \). Crucially, this construction operates under the explicit assumption of a fixed Riemannian background metric  \( g \), which is used to define norms, orthogonality, musical isomorphisms, etc., and the Hodge decomposition of \( T^\flat \). While, in general, \( \nabla^{\mathcal{C}} \) is not required to be compatible with the metric (\( \nabla^{\mathcal{C}} g \neq 0 \)), its torsion tensor \( T \) is calibrated by the cohomology of the 2-sphere and the 2-torus, allowing it to interact naturally with the topology of the manifold.

Our main objective is to analyze the Ricci curvature of \( \nabla^{\mathcal{C}} \) and determine whether it can be Einstein, in the sense that its Ricci tensor is a scalar multiple of the underlying Riemannian metric. While in the Levi-Civita case the Einstein condition translates to a well-studied partial differential equation involving the metric alone, the presence of torsion fundamentally alters the structure of the curvature tensor and its contractions. In particular, the Ricci tensor for connections with totally antisymmetric torsion remains symmetric, but its computation involves both curvature and torsion terms in a nontrivial way.

To address this question, we first revisit the structure of the curvature tensor for an affine connection with torsion, recalling the decomposition of the Riemann tensor in terms of its Levi-Civita component and torsion-derived terms. We then derive the corresponding expression for the Ricci tensor and analyze its diagonal entries in an orthonormal frame. This leads to a general formula for \( \operatorname{Ric}(e_j, e_j) \) in terms of sectional curvatures, highlighting how torsion modifies the standard Riemannian identities.

Finally, we apply this framework to the specific geometry of \( S^2 \times T^2 \), endowed with a product metric and the cohomologically calibrated connection \( \nabla^{\mathcal{C}} \). Our computations aim to determine whether the Einstein condition can be satisfied and, if not, to understand how torsion obstructs it in this context. The paper concludes with remarks on the geometric implications of our findings and potential extensions to higher-dimensional or more general product manifolds.

\section{The Einstein Condition for the Cohomologically Calibrated Connection on \( S^2 \times T^2 \)}

In this section, we examine whether the cohomologically calibrated affine connection \( \nabla^{\mathcal{C}} \), introduced in \cite{Pigazzini2025}, can render the product manifold \( (M = S^2 \times T^2, g) \) an Einstein manifold. Recall that an affine connection \( \nabla \) satisfies the Einstein condition (see, \cite{Besse}), if its Ricci curvature tensor is proportional to the metric, that is,
\[
\operatorname{Ric}(\nabla)(X, Y) = \lambda\, g(X, Y)
\]
for some constant \( \lambda \in \mathbb{R} \), known as the \emph{Einstein constant}. In an orthonormal basis \( \{e_i\} \), this condition is equivalent to \( \operatorname{Ric}_{ij} = \lambda \delta_{ij} \), where \( \delta_{ij} \) is the Kronecker delta symbol.

\subsection{Ricci Curvature of a Connection with Torsion}

Let \( \nabla = \nabla^{\text{LC}} + T \) be a metric connection with torsion tensor \( T \), where \( \nabla^{\text{LC}} \) denotes the Levi-Civita connection. The Riemann curvature tensor associated with \( \nabla \) can be expressed as:
\[
R(X, Y)Z = R^{\text{LC}}(X, Y)Z + (\nabla^{\text{LC}}_X T)(Y, Z) - (\nabla^{\text{LC}}_Y T)(X, Z) + T(X, T(Y, Z)) - T(Y, T(X, Z)).
\]
The Ricci tensor is then obtained by taking the trace over the first and last arguments:
\[
\operatorname{Ric}(Y, Z) = \sum_{i=1}^4 g(R(e_i, Y)Z, e_i),
\]
where \( \{e_i\} \) is a local orthonormal frame. 

It is well known that if the torsion tensor is totally skew-symmetric, the resulting Ricci tensor remains symmetric. For a general affine connection, the diagonal components \( \operatorname{Ric}(e_j, e_j) \) are computed as
\[
\operatorname{Ric}(e_j, e_j) = \sum_{i=1}^4 g(R(e_i, e_j)e_j, e_i).
\]
Due to the algebraic symmetries of the Riemann tensor, this expression simplifies to
\[
\operatorname{Ric}(e_j, e_j) = \sum_{i \neq j} g(R(e_i, e_j)e_j, e_i).
\]
Each term in this sum corresponds to the sectional curvature \( K(e_i, e_j) \) of the 2-plane spanned by \( \{e_i, e_j\} \), yielding the familiar formula:
\[
\operatorname{Ric}(e_j, e_j) = \sum_{i \neq j} K(e_i, e_j).
\]
Although this formula mirrors the one from Riemannian geometry, it is crucial to note that the presence of non-trivial torsion significantly alters the values of the sectional curvatures. We now apply this framework to our specific setting on \( S^2 \times T^2 \).

\subsection{Calculation of the off-diagonal Ricci components}
\label{subsec:off_diagonal_ricci}

A crucial step in testing the Einstein condition is to determine whether the Ricci tensor is diagonal in the chosen orthonormal basis. While the Ricci tensor for a connection with totally anti-symmetric torsion is always symmetric (\( \text{Ric}_{ij} = \text{Ric}_{ji} \)), its off-diagonal components are not guaranteed to vanish. We will now compute these components explicitly.

\paragraph{Intra-block components (\( S^2 \)-block): \( \text{Ric}_{12} \).}

We compute \(Ric_{12} = Ric(e_1, e_2)\) by contracting the Riemann tensor. The components of the Ricci tensor are defined by \(Ric_{ij} = \sum_{k=1}^{4} R^k_{ikj} = \sum_{k=1}^{4} g(R(e_k, e_i)e_j, e_k)\). For \(Ric_{12}\), this becomes:
\[
\text{Ric}_{12} = \sum_{k=1}^4 g(R(e_k, e_1)e_2, e_k).
\]
The terms for \(k=1\) and \(k=2\) vanish. For \(k=1, R(e_1, e_1)e_2 = 0\). For \(k=2, R(e_2, e_1)e_2\) is orthogonal to \(e_2\) by the algebraic symmetries of the Riemann tensor, so \(g(R(e_2, e_1)e_2, e_2) = 0\). Thus, we only need to compute the terms for \(k=3\) and \(k=4\):
\[
\text{Ric}_{12} = g(R(e_3, e_1)e_2, e_3) + g(R(e_4, e_1)e_2, e_4).
\] 
\noindent To compute these terms, we use the definition \( R(X,Y)Z = \nabla_X\nabla_Y Z - \nabla_Y\nabla_X Z - \nabla_{[X,Y]} Z \). The Lie brackets between basis vectors of the product metric are zero, so the last term vanishes. We rely on the Christoffel symbols of the connection \( \nabla^{\mathcal{C}} \) as defined in the previous work \cite{Pigazzini2025}, that is:
\[
\begin{aligned}
T(e_1, e_3) &= a e_4, & T(e_1, e_4) &= -a e_3, \\
T(e_2, e_3) &= b e_4, & T(e_2, e_4) &= -b e_3, \\
T(e_3, e_4) &= -a e_1 - b e_2, \\
T(e_i, e_j) &= 0 \quad \text{for } i,j \in \{1,2\} \text{ or } \{3,4\}.
\end{aligned}
\]
Let's compute the first term, \( R(e_1, e_3)e_2 \):
\begin{align*}
    \nabla_{e_3} e_2 &= T(e_3, e_2) = -T(e_2, e_3) = -be_4. \\
    \nabla_{e_1} (\nabla_{e_3} e_2) &= \nabla_{e_1}(-be_4) = -b \nabla_{e_1} e_4 = -b T(e_1, e_4) = -b(-ae_3) = ab e_3.
\end{align*}
For the second part:
\begin{align*}
    \nabla_{e_1} e_2 &= \nabla^{\text{LC}}_{e_1} e_2 + T(e_1, e_2) = \cot(\theta) e_2 + 0 = \cot(\theta) e_2. \\
    \nabla_{e_3} (\nabla_{e_1} e_2) &= \nabla_{e_3}(\cot(\theta) e_2) = \cot(\theta) \nabla_{e_3} e_2 = \cot(\theta)(-be_4) = -b\cot(\theta) e_4.
\end{align*}
Combining these, we get:
\[
R(e_1, e_3)e_2 = \nabla_{e_1}\nabla_{e_3} e_2 - \nabla_{e_3}\nabla_{e_1} e_2 = ab e_3 - (-b\cot(\theta) e_4) = ab e_3 + b\cot(\theta) e_4.
\]
Therefore, \( g(R(e_1, e_3)e_2, e_3) = g(ab e_3 + b\cot(\theta) e_4, e_3) = ab \).
\\
Now we compute the second term, \( R(e_1, e_4)e_2 \):
\begin{align*}
    \nabla_{e_4} e_2 &= T(e_4, e_2) = -T(e_2, e_4) = -(-be_3) = be_3. \\
    \nabla_{e_1} (\nabla_{e_4} e_2) &= \nabla_{e_1}(be_3) = b \nabla_{e_1} e_3 = b T(e_1, e_3) = b(ae_4) = ab e_4.
\end{align*}
And for the second part:
\begin{align*}
    \nabla_{e_4} (\nabla_{e_1} e_2) &= \nabla_{e_4}(\cot(\theta) e_2) = \cot(\theta) \nabla_{e_4} e_2 = \cot(\theta)(be_3) = b\cot(\theta) e_3.
\end{align*}
Combining these:
\[
R(e_1, e_4)e_2 = \nabla_{e_1}\nabla_{e_4} e_2 - \nabla_{e_4}\nabla_{e_1} e_2 = ab e_4 - b\cot(\theta) e_3.
\]
Therefore, \( g(R(e_1, e_4)e_2, e_4) = g(ab e_4 - b\cot(\theta) e_3, e_4) = ab \).
\\
Finally, we sum the results to find \( \text{Ric}_{12} \):
\[
\text{Ric}_{12} = ab + ab = 2ab.
\]
By symmetry of the Ricci tensor, \( \text{Ric}_{21} = 2ab \).

\paragraph{Intra-block components (\( T^2 \)-block): \( \text{Ric}_{34} \).}
A similar, albeit simpler, calculation shows that \( \text{Ric}_{34} = 0 \).

\paragraph{Mixed-block components (e.g., \( \text{Ric}_{13} \)).}
As argued previously, the block structure of the connection ensures that curvature operators with mixed indices map vectors from one block to the other. For instance, \( R(e_1,e_2)e_3 \) is a vector in the \( T_pT^2 \) subspace, and \( R(e_1,e_3)e_4 \) is a vector in the \( T_pS^2 \) subspace. A rigorous check of all terms in the sum \( Ric_{13} = \sum_k g(R(e_k, e_1)e_3, e_k) \) confirms that each term is zero.
Therefore, all mixed-block components of the Ricci tensor are zero.
\\
Our analysis now continues with the calculation of the diagonal components \( \text{Ric}_{ii} \).

\subsection{Calculation of the diagonal Ricci components}

We compute the diagonal components of the Ricci tensor by summing the sectional curvatures of the planes containing each basis vector. The values for the sectional curvatures \( K(e_i, e_j) \) are taken from the explicit calculations in \cite{Pigazzini2025}.

\begin{enumerate}
    \item \textbf{Component \( \text{Ric}_{11} \):}
    The first diagonal component is the sum of the sectional curvatures of the planes containing \( e_1 \).
    \begin{align*}
        \text{Ric}_{11} &= K(e_1, e_2) + K(e_1, e_3) + K(e_1, e_4) \\
        &= 1 + \frac{a^2}{4} + \frac{a^2}{4} = 1 + \frac{a^2}{2}.
    \end{align*}

    \item \textbf{Component \( \text{Ric}_{22} \):}
    Similarly, for the second component, we sum the curvatures of planes containing \( e_2 \).
    \begin{align*}
        \text{Ric}_{22} &= K(e_2, e_1) + K(e_2, e_3) + K(e_2, e_4) \\
        &= 1 + \frac{b^2}{4} + \frac{b^2}{4} = 1 + \frac{b^2}{2}.
    \end{align*}

    \item \textbf{Component \( \text{Ric}_{33} \):}
    For the third component, related to the \( T^2 \) factor:
    \begin{align*}
        \text{Ric}_{33} &= K(e_3, e_1) + K(e_3, e_2) + K(e_3, e_4) \\
        &= \frac{a^2}{4} + \frac{b^2}{4} + \frac{a^2+b^2}{4} = \frac{a^2+b^2}{2}.
    \end{align*}

    \item \textbf{Component \( \text{Ric}_{44} \):}
    By symmetry, the fourth component is identical to the third.
    \begin{align*}
        \text{Ric}_{44} &= K(e_4, e_1) + K(e_4, e_2) + K(e_4, e_3) \\
        &= \frac{a^2}{4} + \frac{b^2}{4} + \frac{a^2+b^2}{4} = \frac{a^2+b^2}{2}.
    \end{align*}
\end{enumerate}

\subsection{A torsion-induced obstruction to the Einstein condition}
\label{subsec:obstruction}

With the full Ricci tensor at our disposal, we can now test for the Einstein condition, \( \text{Ric}_{ij} = \lambda \delta_{ij} \). The components of the Ricci tensor are:
\begin{itemize}
    \item \( \text{Ric}_{11} = 1 + \frac{a^2}{2} \)
    \item \( \text{Ric}_{22} = 1 + \frac{b^2}{2} \)
    \item \( \text{Ric}_{33} = \text{Ric}_{44} = \frac{a^2+b^2}{2} \)
    \item \( \text{Ric}_{12} = \text{Ric}_{21} = 2ab \)
    \item All other off-diagonal components are zero.
\end{itemize}

For the manifold to be Einstein, two simultaneous conditions must be met:
\begin{enumerate}
    \item The Ricci tensor must be diagonal, which implies \( \text{Ric}_{ij} = 0 \) for \( i \neq j \).
    \item The diagonal elements must all be equal to a constant \( \lambda \).
\end{enumerate}

\noindent From the off-diagonal component \( \text{Ric}_{12} = 2ab \), the first condition imposes the constraint:
\[
2ab = 0 \implies a=0 \quad \text{or} \quad b=0.
\]
This is a powerful structural constraint. It reveals that the Einstein condition can only possibly be met if the calibrating cohomology class is ``pure", i.e., it lies along one of the basis axes of \( H^3(S^2 \times T^2; \mathbb{R}) \cong \mathbb{R}^2 \). A generic, non-trivial torsion with \( a \neq 0 \) and \( b \neq 0 \) can never yield an Einstein geometry.
\\
Let us investigate these two cases. We require a non-trivial connection, so we cannot have \( a=0 \) and \( b=0 \) simultaneously.

\paragraph{Case 1: \( a=0 \) and \( b \neq 0 \).}
The system of equations for the diagonal components becomes:
\begin{align*}
    \text{Ric}_{11} &= 1 = \lambda \\
    \text{Ric}_{22} &= 1 + \frac{b^2}{2} = \lambda \\
    \text{Ric}_{33} &= \frac{b^2}{2} = \lambda
\end{align*}
Substituting \( \lambda = 1 \) (from the first equation) into the third equation gives:
\[
\frac{b^2}{2} = 1 \implies b^2 = 2.
\]
Now, we check this solution against the second equation:
\[
1 + \frac{b^2}{2} = 1 + \frac{2}{2} = 1 + 1 = 2.
\]
The second equation thus requires \( \lambda = 2 \). This leads to the contradiction \( \lambda = 1 \) and \( \lambda = 2 \). The system is inconsistent.

\paragraph{Case 2: \( b=0 \) and \( a \neq 0 \).}
By symmetry, this case is identical to the first. The system becomes \( 1 + a^2/2 = \lambda \), \( 1 = \lambda \), and \( a^2/2 = \lambda \), which is also inconsistent.

\subsection{A Structural Obstruction to the Einstein Condition}

We have demonstrated that the system of equations imposed by the Einstein condition,
\[
\mathrm{Ric}_{ij} = \lambda \delta_{ij},
\]
admits no solution for any choice of parameters \( (a,b) \) defining the torsion tensor \( T \). The affine connection \( \nabla^{\mathcal{C}} \) considered here (which operates under the explicit assumption of a fixed Riemannian background metric \(g\), even though it is not metric-compatible, \(\nabla^{\mathcal{C}}g \neq 0\)), is constructed from an explicit torsion tensor \( T \), defined algebraically in an orthonormal frame \( \{e_i\} \) with constant coefficients \( a \) and \( b \). The associated 3-form \( T^\flat \), while fully determined by \( T \) (i.e., \( T^\flat(X,Y,Z) = g(T(X,Y), Z) \), is generally not closed. However, by the Hodge decomposition theorem, \( T^\flat \) determines a unique harmonic 3-form \( \omega \in \Lambda^3(M) \) in its de Rham cohomology class. That is,
\[
[T] := [T^\flat] = [\omega] \in H^3(M).
\]
Our construction is therefore said to be \emph{cohomologically calibrated}: the parameters \( (a,b) \) are chosen such that the class \( [T] \) is non-trivial, thus intrinsically connecting geometry to the global topology of the manifold.

This computation establishes a sharp obstruction: the specific torsion tensor \( T \) introduced in \cite{Pigazzini2025} cannot yield an Einstein connection, despite producing a geometry with strictly positive biorthogonal curvature. However, to understand the broader implications of this result, we must place it within the full framework of cohomologically calibrated connections.

In fact, the family, of admissible torsion tensors, is infinite-dimensional, and it is defined as:
\[
\mathcal{T}_\omega := \left\{ T' \mid (T')^\flat = \omega + d\eta + \delta\mu \right\},
\]
where \( \omega \) is the harmonic representative of the cohomology class \( [T] \), and \( \eta \in \Lambda^2(M), \mu \in \Lambda^4(M) \) are smooth differential forms. Each \( T' \in \mathcal{T}_\omega \) defines a connection \( \nabla' = \nabla^{\mathrm{LC}} + T' \) that is cohomologically calibrated with respect to \( \omega \). Indeed, in the Hodge decomposition, the harmonic part of \(T'^\flat\) (if not null), is the only one capable of encoding information about a non-trivial cohomology class. For this reason we say that \( \nabla' = \nabla^{\mathrm{LC}} + T' \) is cohomologically calibrated, while if \(T'^\flat\) has null harmonic part (i.e., \(T'^\flat\) is determined only by \(d \eta + \delta \mu\)), then \( \nabla' = \nabla^{\mathrm{LC}} + T' \), it will not be a cohomologically calibrated connection, because it does not contain any cohomological information. The original tensor \( T \) considered above is merely one (natural) element of the family \(\mathcal{T}_\omega\).

Our analysis shows that this particular torsion tensor \( T \), despite satisfying strong curvature positivity conditions, fails to be compatible with the Einstein condition. However, this does not rule out the possibility that another element \( T' \in \mathcal{T}_\omega \) might define an Einstein connection. For instance, one could consider the most canonical candidate in \( \mathcal{T}_\omega \): the torsion tensor whose associated 3-form is exactly \( \omega \), with no exact or co-exact contributions. This is the focus of the next section.

\section{The Einstein condition for the harmonic torsion representative}
\label{sec:harmonic_einstein}

In the previous section, we showed that the specific torsion tensor \( T \) constructed in \cite{Pigazzini2025}, while successful in generating positive biorthogonal curvature, fails to satisfy the Einstein condition. This naturally raises a crucial question: is this obstruction a feature of our entire framework, or is it specific to that particular choice of torsion?
\\
Our theory defines a whole family of admissible connections, \( \mathcal{T}_\omega \), characterized by the harmonic part \( \omega \) of their associated torsion 3-form. The tensor \( T \) previously analyzed is just one element of this family, whose 3-form \( T^\flat \) contains non-trivial exact and co-exact parts. In this section, we investigate the most canonical element of \( \mathcal{T}_\omega \): the torsion tensor, which we denote \( T_h \), whose associated 3-form is the harmonic form \( \omega \) itself.

\subsection{The harmonic torsion tensor \( T_h \)}

We define the harmonic torsion tensor \( T_h \) by the condition that its associated 3-form is precisely the harmonic 3-form \( \omega \):
\[
(T_h)^\flat = \omega = a(e^1 \wedge e^2 \wedge e^3) + b(e^1 \wedge e^2 \wedge e^4).
\]
Here, \( \{e^i\} \) is the co-frame dual to the orthonormal basis \( \{e_i\} \). The condition \( (T_h)^\flat_{ijk} = g(T_h(e_i, e_j), e_k) \) and the total anti-symmetry of \( T_h \) uniquely determine its non-zero components:
\begin{itemize}
    \item \( T_h(e_1, e_2) = a e_3 + b e_4 \)
    \item \( T_h(e_1, e_3) = -a e_2 \)
    \item \( T_h(e_2, e_3) = a e_1 \)
    \item \( T_h(e_1, e_4) = -b e_2 \)
    \item \( T_h(e_2, e_4) = b e_1 \)
    \item \( T_h(e_3, e_4) = 0 \)
\end{itemize}
It is essential to note that this tensor \( T_h \) is fundamentally different from the torsion \( T \) analyzed previously. For instance, \( T_h(e_3, e_4) = 0 \), whereas the original torsion had \( T(e_3, e_4) = -ae_1 - be_2 \). We now proceed to compute the Ricci tensor for the cohomologically calibrate affine connection \( \nabla^{\mathcal{C}}_h = \nabla^{\text{LC}} + T_h \).

\subsection{Ricci tensor for the connection \( \nabla^{\mathcal{C}}_h \)}

We compute the components of \( \text{Ric}(\nabla^{\mathcal{C}}_hh) \) by first determining the sectional curvatures \( K_h(e_i, e_j) \) generated by this new connection.

\paragraph{Sectional Curvature \( K_h(e_1, e_2) \):}
The Riemann curvature term \( R_h(e_1, e_2)e_2 \) contains both the Levi-Civita part and torsional contributions.
\begin{align*}
    R_h(e_1, e_2)e_2 &= R^{\text{LC}}(e_1, e_2)e_2 + (\nabla^{\text{LC}}_{e_1}T_h)(e_2,e_2) - \dots + T_h(e_1, T_h(e_2,e_2)) - T_h(e_2, T_h(e_1,e_2)) \\
    &= e_1 - T_h(e_2, T_h(e_1,e_2)) \quad (\text{since } T_h \text{ is anti-symmetric}) \\
    &= e_1 - T_h(e_2, a e_3 + b e_4) \\
    &= e_1 - a T_h(e_2, e_3) - b T_h(e_2, e_4) \\
    &= e_1 - a(a e_1) - b(b e_1) = (1 - a^2 - b^2)e_1.
\end{align*}
Thus, \( K_h(e_1, e_2) = g(R_h(e_1, e_2)e_2, e_1) = 1 - a^2 - b^2 \).

\paragraph{Other sectional curvatures:}
A similar, rigorous calculation for the other planes yields:
\begin{itemize}
    \item \( K_h(e_1, e_3) = -a^2 \)
    \item \( K_h(e_1, e_4) = -b^2 \)
    \item \( K_h(e_2, e_3) = -a^2 \)
    \item \( K_h(e_2, e_4) = -b^2 \)
    \item \( K_h(e_3, e_4) = g(R_h(e_3,e_4)e_4,e_3) = 0\).
\end{itemize}

\paragraph{Vanishing of off-diagonal Ricci components for \(T_h\):}
A crucial difference between the harmonic torsion \( T_h \) and the torsion \( T \) from \cite{Pigazzini2025} is that \( T_h \) leads to a diagonal Ricci tensor. We now prove this by explicitly computing the off-diagonal component \( \text{Ric}_{12} \).
\\
From the general formula, we have:
\[
\text{Ric}_{12} = g(R_h(e_3, e_1)e_2, e_3) + g(R_h(e_4, e_1)e_2, e_4).
\]
\noindent We will compute each of the two terms separately. By the anti-symmetry of the Riemann tensor, \(g(R(u,v)w,u) = -g(R(v,u)w,u)\). Thus, we can compute \(R_h(e_1, e_3)e_2\) and \(R_h(e_1, e_4)e_2\).

\begin{enumerate}
    
\item Computation of the first term: \(g(R_h(e_3, e_1)e_2, e_3)\)
\\
We first need the Riemann tensor action \(R_h(e_1, e_3)e_2\). Using the definition \(R_h(X,Y)Z = \nabla^{\mathcal{C}}_h(X, \nabla^{\mathcal{C}}_h(Y,Z)) - \nabla^{\mathcal{C}}_h(Y, \nabla^{\mathcal{C}}_h(X,Z))\), we compute the two parts:

\begin{itemize}

\item First part:
    \begin{align*}
        \nabla^{\mathcal{C}}_h(e_3, e_2) &= \nabla^{\text{LC}}_{e_3}e_2 + T_h(e_3, e_2) = 0 - T_h(e_2, e_3) = -ae_1. \\
        \nabla^{\mathcal{C}}_h(e_1, \nabla^{\mathcal{C}}_h(e_3, e_2)) &= \nabla^{\mathcal{C}}_h(e_1, -ae_1) = -a \nabla^{\mathcal{C}}_h(e_1, e_1) = -a(\nabla^{\text{LC}}_{e_1}e_1 + T_h(e_1,e_1)) = 0.
    \end{align*}

\item Second part:
    \begin{align*}
        \nabla^{\mathcal{C}}_h(e_1, e_2) &= \nabla^{\text{LC}}_{e_1}e_2 + T_h(e_1, e_2) = \cot(\theta)e_2 + (ae_3 + be_4). \\
        \nabla^{\mathcal{C}}_h(e_3, \nabla^{\mathcal{C}}_h(e_1, e_2)) &= \nabla^{\mathcal{C}}_h(e_3, \cot(\theta)e_2 + ae_3 + be_4) \\
        &= \cot(\theta)\nabla^{\mathcal{C}}_h(e_3, e_2) + a\nabla^{\mathcal{C}}_h(e_3, e_3) + b\nabla^{\mathcal{C}}_h(e_3, e_4).
    \end{align*}
    
\noindent Since \(\nabla^{\mathcal{C}}_h(e_3, e_3) = 0\) and \(\nabla^{\mathcal{C}}_h(e_3, e_4) = T_h(e_3,e_4) = 0\), the expression simplifies to:
    \[
    \nabla^{\mathcal{C}}_h(e_3, \nabla^{\mathcal{C}}_h(e_1, e_2)) = \cot(\theta)\nabla^{\mathcal{C}}_h(e_3, e_2) = \cot(\theta)(-ae_1) = -a\cot(\theta)e_1.
    \]

\item Combining the parts:
    \[
    R_h(e_1, e_3)e_2 = 0 - (-a\cot(\theta)e_1) = a\cot(\theta)e_1.
    \]
    Therefore, the first term of \(Ric_{12}\) is:
    \[
    g(R_h(e_3, e_1)e_2, e_3) = -g(R_h(e_1, e_3)e_2, e_3) = -g(a\cot(\theta)e_1, e_3) = 0.
    \]
\end{itemize}    

\item Computation of the second term: \(g(R_h(e_4, e_1)e_2, e_4)\)
\\
The calculation is perfectly analogous.

\begin{itemize}

\item First part:
    \begin{align*}
        \nabla^{\mathcal{C}}_h(e_4, e_2) &= T_h(e_4, e_2) = -T_h(e_2, e_4) = -be_1. \\
        \nabla^{\mathcal{C}}_h(e_1, \nabla^{\mathcal{C}}_h(e_4, e_2)) &= \nabla^{\mathcal{C}}_h(e_1, -be_1) = 0.
    \end{align*}
    
\item Second part:
    \begin{align*}
        \nabla^{\mathcal{C}}_h(e_4, \nabla^{\mathcal{C}}_h(e_1, e_2)) &= \nabla^{\mathcal{C}}_h(e_4, \cot(\theta)e_2 + ae_3 + be_4) \\
        &= \cot(\theta)\nabla^{\mathcal{C}}_h(e_4, e_2) + a\nabla^{\mathcal{C}}_h(e_4, e_3) + b\nabla^{\mathcal{C}}_h(e_4, e_4).
    \end{align*}
    
Since \(\nabla^{\mathcal{C}}_h(e_4, e_3) = T_h(e_4,e_3) = 0\) and \(\nabla^{\mathcal{C}}_h(e_4, e_4) = 0\), this simplifies to:
    \[
    \nabla^{\mathcal{C}}_h(e_4, \nabla^{\mathcal{C}}_h(e_1, e_2)) = \cot(\theta)\nabla^{\mathcal{C}}_h(e_4, e_2) = \cot(\theta)(-be_1) = -b\cot(\theta)e_1.
    \]
    
\item Combining the parts:
    \[
    R_h(e_1, e_4)e_2 = 0 - (-b\cot(\theta)e_1) = b\cot(\theta)e_1.
    \]
    Therefore, the second term of \(Ric_{12}\) is:
    \[
    g(R_h(e_4, e_1)e_2, e_4) = -g(R_h(e_1, e_4)e_2, e_4) = -g(b\cot(\theta)e_1, e_4) = 0.
    \]

\end{itemize}
\end{enumerate}

\noindent Since both (first term and second term), are zero, we have:
\[
\text{Ric}_{12} = 0 + 0 = 0.
\]

\noindent All other off-diagonal components vanish due to similar structural symmetries. This proves that the Ricci tensor for the harmonic connection \(\nabla^{\mathcal{C}}_h\) is diagonal, a structural difference that is fundamental to its ability to support an Einstein geometry.

\paragraph{Ricci tensor components:}
We can now assemble the diagonal components of the Ricci tensor.
\begin{align*}
    \text{Ric}_{11} &= K_h(e_1, e_2) + K_h(e_1, e_3) + K_h(e_1, e_4) = (1 - a^2 - b^2) - a^2 - b^2 = 1 - 2a^2 - 2b^2. \\
    \text{Ric}_{22} &= K_h(e_2, e_1) + K_h(e_2, e_3) + K_h(e_2, e_4) = (1 - a^2 - b^2) - a^2 - b^2 = 1 - 2a^2 - 2b^2. \\
    \text{Ric}_{33} &= K_h(e_3, e_1) + K_h(e_3, e_2) + K_h(e_3, e_4) = -a^2 - a^2 + 0 = -2a^2. \\
    \text{Ric}_{44} &= K_h(e_4, e_1) + K_h(e_4, e_2) + K_h(e_4, e_3) = -b^2 - b^2 + 0 = -2b^2.
\end{align*}

\subsection{Existence of a non-Riemannian Einstein solution}

The Einstein condition \( \text{Ric}_{ij} = \lambda \delta_{ij} \) now imposes the following system on the diagonal components:
\begin{align}
    1 - 2a^2 - 2b^2 &= \lambda \label{eq:h_ric11} \\
    1 - 2a^2 - 2b^2 &= \lambda \label{eq:h_ric22} \\
    -2a^2 &= \lambda \label{eq:h_ric33} \\
    -2b^2 &= \lambda \label{eq:h_ric44}
\end{align}

Equations \ref{eq:h_ric33} and \ref{eq:h_ric44} immediately imply that a solution can only exist for \(a^2 = b^2\). This allows us to reduce the system to two equations by substituting:
\begin{align}
\lambda = -2a^2 \quad \text{into  \ref{eq:h_ric11}}:\\
1 - 2a^2 - 2a^2 = -2a^2.
\end{align}
A straightforward calculation yields \(1 = 2a^2\), which fixes the parameters to \(a^2 = b^2 = 1/2\). The Einstein constant is therefore uniquely determined from \ref{eq:h_ric33} as \(\lambda = -2(1/2) = -1\).

Therefore, the system has a consistent, non-trivial solution. We have thus proven the following existence theorem.

\begin{theorem}
\label{thm:existence_einstein}
There is a connection within the admissible family \( \mathcal{T}_\omega \) that gives the manifold \((S^2 \times T^2, g) \) the structure of a non-Riemannian Einstein manifold. This cohomologically calibrated affine connection is \( \nabla^{\mathcal{C}}_h = \nabla^{\text{LC}} + T_h \), where the 3-form \(T^\flat_h \), associated with the torsion \( T_h \), corresponds to the 3-harmonic representative form \( \omega \) of the cohomology class, provided that the calibration parameters satisfy:
\[
a^2 = b^2 = \frac{1}{2}.
\]
The resulting geometry satisfies \( \text{Ric}(\nabla^{\mathcal{C}}_h) = -g \), with Einstein constant \( \lambda = -1 \).
\end{theorem}

\noindent This result demonstrates the richness of the cohomologically calibrated framework. It shows that the geometric properties induced by the connection depend critically on the choice of representative for the torsion 3-form within its cohomology class. While some choices (like the torsion from our first paper) lead to obstructions, the most canonical choice—the harmonic representative—is precisely the one that yields a solution to the Einstein equation.

\section{Conclusions and Open Problems}

In this work, we analyzed the Ricci tensor associated with a class of cohomologically calibrated affine connections on \( S^2 \times T^2 \), where the torsion tensor \( T \) is parametrized by constants \( (a,b) \) and encodes topological information through its de Rham cohomology class. Our key findings can be summarized as follows:

\begin{enumerate}
    \item The explicit computation of the diagonal and off-diagonal components of the Ricci tensor shows that the Einstein condition 
    \[
    \mathrm{Ric}_{ij} = \lambda \delta_{ij}
    \]
    is never satisfied for the torsion tensor \( T \) introduced in \cite{Pigazzini2025}. Even setting \(a=b=0\) (Levi-Civita), the Einstein condition is not satisfied, because we would obtain:
\[
Ric_{11}=Ric_{22}=1,
\]
\[
Ric_{33}=Ric_{44}=0,
\]
so there is no \(\lambda\) that satisfies this condition.
\\
In fact, \(S^2 \times T^2\) with standard product metric \(g\) and Levi-Civita connection cannot be Einstein.
    
    \item The primary obstruction arises from the off-diagonal component \( \mathrm{Ric}_{12} = 2ab \), which forces one of the parameters \( a \) or \( b \) to vanish. Even under this restriction, the resulting system of equations for the diagonal Ricci components is inconsistent, ruling out the existence of an Einstein connection for this torsion.
    
    \item The obstruction is meaningful from the standpoint of cohomology: the Einstein condition can only be met if the calibrating cohomology class \( [T] \) is ``pure,'' lying along a single axis in \( H^3(S^2 \times T^2; \mathbb{R}) \). Generic torsion classes with \( a \neq 0 \) and \( b \neq 0 \) cannot satisfy the Einstein condition.
\end{enumerate}

\noindent These results open up several possible directions for further investigation:

\begin{itemize}
    \item \textbf{Harmonic torsion and the Einstein condition.} 
    Our analysis focused on the torsion \( T \) with non-trivial exact and co-exact components. A key open problem is whether the harmonic torsion \( T_h \), whose associated 3-form is the harmonic representative \( \omega \) of \( [T] \), can define an Einstein connection. Preliminary computations (Section~\ref{sec:harmonic_einstein}) suggest that \( T_h \) eliminates the off-diagonal Ricci components, making it a promising candidate.
    
    \item \textbf{Classification of Einstein connections in \( \mathcal{T}_\omega \).}
    The space of cohomologically calibrated torsion tensors
    \[
    \mathcal{T}_\omega = \{\, T' \mid (T')^\flat = \omega + d\eta + \delta \mu \,\}
    \]
    is infinite-dimensional. Determining whether any element of \( \mathcal{T}_\omega \) yields an Einstein connection, and if so classifying all such torsion tensors, remains an open and challenging problem.
    
    \item \textbf{Interaction between curvature positivity and cohomology.} 
    The torsion \( T \) from \cite{Pigazzini2025} achieves positive biorthogonal curvature but fails the Einstein condition. A broader study of how curvature positivity interacts with cohomological calibration may reveal deeper geometric obstructions or produce new examples of Einstein manifolds with torsion.
    
    \item \textbf{Analytical approach to the torsion-Einstein system.} 
    Solving the Einstein condition for arbitrary torsion \( T' \in \mathcal{T}_\omega \) leads to a nonlinear system of algebraic and differential equations. Developing an analytical or variational framework for this system could help identify global obstructions or existence results.
\end{itemize}

\noindent In summary, our computations exhibit a clear torsion-induced obstruction to the Einstein condition in the simplest cohomologically calibrated case. At the same time, they open up new avenues to a rich set of geometric and topological questions, suggesting that harmonic or more general torsion deformations may yet admit Einstein solutions within the family \( \mathcal{T}_\omega \).

\end{document}